\documentclass{ifacconf}

\usepackage{amsmath} 
\usepackage{amssymb}  
\usepackage{url}
\usepackage{color}

\newtheorem{proposition}[thm]{Proposition}

\newtheorem{remark}{Remark}
\newtheorem{definition}{Definition}
\newtheorem{assumption}{Assumption}

\usepackage{graphicx}      
\usepackage{natbib}        
\begin{document}
\begin{frontmatter}

\title{Which notion of energy for bilinear quantum systems?\thanksref{footnoteinfo}} 

\thanks[footnoteinfo]{This work has been partially supported by INRIA Nancy-Grand Est.
Second and third authors were partially supported by
French Agence National de la Recherche ANR ``GCM'' program
``BLANC-CSD'', contract number NT09-504590. The
third author was partially supported by European Research
Council ERC StG 2009 ``GeCoMethods'', contract number
239748. }

\author[First]{Nabile Boussa\"{i}d} 
\author[Second]{Marco Caponigro} 
\author[Third]{Thomas Chambrion}

\address[First]{Laboratoire de math\'ematiques,
Universit\'e de Franche--Comt\'e,
25030~Besan\c{c}on, France
(email:  Nabile.Boussaid@univ-fcomte.fr).}
\address[Second]{
Dept of Mathematical Sciences and
 Center for Computational and Integrative Biology,
Rutgers University,
08102~Camden, NJ, USA
 (email: marco.caponigro@rutgers.edu).}
\address[Third]{Universit\'e de Lorraine, Institut \'Elie Cartan de Nancy, UMR 7502,  Vand{\oe}uvre,  F-54506, France and
INRIA, Villers-l\`es-Nancy,F-54600, France
(email: Thomas.Chambrion@inria.fr).}

\begin{abstract}                
In this note we investigate 
what is the best $L^p$-norm in order to describe the relation between the evolution of the state of a bilinear quantum system with the $L^p$-norm of the external field. Although $L^{2}$ has a structure more easy to handle, the $L^{1}$ norm is more suitable for this purpose.
Indeed for every $p>1$ it is possible to steer, with arbitrary precision, a generic bilinear quantum system from any eigenstate of the free Hamiltonian to any other with a control of arbitrary small $L^p$ norm. Explicit optimal costs for the $L^1$ norm are computed on an example. 
\end{abstract}


\end{frontmatter}

\section{INTRODUCTION}
\subsection{Physical context}

The state of a quantum system evolving on a Riemannian
manifold $\Omega$, with associated measure $\mu$, is described by its \emph{wave
function}, that is, a point in the unit sphere of $L^2(\Omega, \mathbf{C})$. A
system with wave function $\psi$ is in a subset $\omega$ of $\Omega$ with
probability $\displaystyle{\int_{\omega} \!\!\!|\psi|^2 \mathrm{d}\mu}$.

 When submitted to an
excitation by an external  field ({\it e.g.} a laser) the time evolution of the wave 
function is governed by the bilinear Schr\"odinger
equation 
\begin{equation}\label{eq:blse}
\mathrm{i} \frac{\partial \psi}{\partial t}=-\Delta \psi +V(x)
\psi(x,t) +u(t) W(x) \psi(x,t),
\end{equation}
where $V, W:\Omega\rightarrow \mathbf{R}$ are real functions describing  respectively the
physical properties of the uncontrolled system and the external field, and 
$u:\mathbf{R}\rightarrow \mathbf{R}$ is a real function of the time representing the intensity of the latter.

\subsection{Energy for a quantum system}

Physically, the energy of a quantum system (\ref{eq:blse}) with wave function $\psi$ is $E(\psi)=\int_{\Omega}\left \lbrack (-\Delta +V) \overline{\psi } \right \rbrack \psi~ \mathrm{d}\mu$. 
The energy is therefore constant in time whenever the control $u$ is zero. When the control $u$ is nonzero, and provided suitable regularity hypotheses, the energy evolves as 
\begin{equation}\label{EQ_evolution_energy_blse}
\frac{\mathrm{d}E}{\mathrm{d}t}=2 u(t) \Im \left ( \int_{\Omega}\left \lbrack  (\Delta+V)\overline{\psi} \right \rbrack W\psi~ \mathrm{d}\mu \right ).
\end{equation}
Note that the time derivative of the energy $E$ at time $t$ depends on the value $u(t)$ of the intensity of the external field  \emph{and} on the wave function $\psi(t)$.
Therefore a natural question is to find an \emph{a priori} relation between the time evolution of the energy of  system (\ref{eq:blse}) and properties of the external field represented by $u$. 
In particular we address the problem of finding a bound on the energy after the action of an external field. Namely, given an initial condition $\psi_0$ and a control $u:[0,T]\rightarrow \mathbf{R}$
denoting with $\psi$ the solution of (\ref{eq:blse}) with $\psi(0)=\psi_0$ we look for 
bounds on the energy $E(T)$ in terms of the $L^p$ norm of $u$
 $$
 \|u\|_{L^p(0,T)}=\left ( \int_0^T |u(t)|^p \mathrm{d}t \right )^{\frac{1}{p}},
 $$
 (for some suitable $p>0$) without computing explicitly the solution $t\in [0,T]\mapsto \psi(t)$.

{Many previous works adressed the problem of the optimal control of the system (\ref{eq:blse}) for costs involving the 
$L^2$ norm of the control 
(see for instance \cite{PhysRevA.42.1065} or \cite{symeon}).
 The main reason for the choice of the $L^2$ norm is the fact that the natural Hilbert structure of $L^2$ spaces 
allows to use the powerful tools of Hilbert optimization.}
{{It is sometimes believed}} that there is a natural relation of the $L^2$ norm of $u$ and the energy of the systems.
This note presents a priori bounds on the $L^{p}$-norm of the control and shows that, in general, 
the $L^1$-norm provides more informations on the evolution of the system{{than}} other $L^{p}$-norms for $p>1$. Most of the material presented below is classical for finite dimensional conservative bilinear systems. The contribution of the present note is to treat in a rigorous and unified way the case of both finite and infinite  dimensional systems. 

\subsection{Content of the paper}

The first part of the paper (Section \ref{SEC_Contr_results}) presents some theoretical tools for bilinear quantum systems. The question of estimation of the energy is reformulated in terms of a problem of optimal control. Some relations between the variation of energy of a quantum system and the $L^1$ norm of the external field are 
given in Section \ref{SEC_C_r}. Finally, some explicit computations are presented on an example in Section  \ref{SEC_example}.
\section{INFINITE DIMENSIONAL QUANTUM SYSTEMS}\label{SEC_Contr_results}
\subsection{Abstract framework}

We reformulate the problem (\ref{eq:blse}) in a more abstract framework.  This
will allow us  to treat examples slightly more general than (\ref{eq:blse}),
for instance, the example in~\cite[Section III.A]{weakly-coupled}. In a separable Hilbert
space $H$
endowed with norm $\| \cdot \|$ and Hilbert product $\langle \cdot, \cdot
\rangle$, we consider the evolution problem
\begin{equation}\label{eq:main}
\frac{d \psi}{dt}=(A+ u(t) B)\psi
\end{equation}
where $(A,B)$ satisfies the following assumption.

\begin{assumption}\label{ass:ass}
$(A,B)$ is a pair of linear operators such that
 \begin{enumerate}
 \item $A$ is skew-adjoint and has purely discrete spectrum $(-\mathrm{i}
\lambda_k)_{k \in \mathbf{N}}$, the sequence  $(\lambda_k)_{k \in \mathbf{N}}$
is positive non-decreasing and accumulates at
$+\infty$;\label{ASS_bounded_from_below}
	\item $B:H\rightarrow H$ is skew-adjoint and bounded.\label{ASS_B_bounded}
 \end{enumerate}
\end{assumption}

In the rest of our study, we denote by $(\phi_k)_{k \in
\mathbf{N}}$ an Hilbert basis of $H$ such
that $A\phi_k=-\mathrm{i}  \lambda_k \phi_k$ for every $k$ in $\mathbf{N}$.
We denote by $D(A+uB)$ the domain where
$A+uB$ is skew-adjoint.

 From Assumption~\ref{ass:ass}.\ref{ASS_B_bounded}
together with Kato-Rellich Theorem,
we deduce that $A+uB$ is skew-adjoint with domain $D(A)$.  Moreover, for every
constant $u$ in $\mathbf{R}$,  $\mathrm{i}(A+ u B )$ is bounded from below.

 Hence, for every initial condition $\psi_0$ in
$H$, 
for every
$u$ piecewise constant, $u:t\mapsto \sum_{j} u_j \chi_{(t_j,t_{j+1})}(t)$, with $0=t_0\leq t_1 \leq \ldots  \leq t_{N+1}$ and $u_0,\ldots,u_N$ in $\mathbf{R}$,
one can define
 the solution $t\mapsto \Upsilon^u_t\psi_0$ of \eqref{eq:main}   by
\begin{multline*}
 \Upsilon^u_t \psi_{0}=e^{(t-t_{j-1})(A + u_{j-1} B)}\circ\\\circ
e^{(t_{j-1}-t_{j-2})(A + u_{j-2} B)}\circ \cdots  \circ e^{t_{0}(A+
u_0  B)} \psi_{0},
\end{multline*}
for $t\in [t_{j-1},t_{j})$.
For a control $u$ in $L^{1}(\mathbf{R})$ we define the solution 
 using the 
 following  continuity result.
 
\begin{proposition}\label{PRO_continuity_Upsilon}
Let $u$ and $(u_n)_{n\in \mathbf{N}}$ be in  $L^1(\mathbf{R})$. If for every $t$ in $\mathbf{R}$ $\int_0^t u_n(\tau) \mathrm{d} \tau$ converges to  $\int_0^t u(\tau) \mathrm{d} \tau$ as $n$ tends to infinity, then, for every $t$ in $\mathbf{R}$ and every $\psi_0$ in $H$, $(\Upsilon^{u_n}_t\psi_0 )_{n \in \mathbf{N}}$ converges to $\Upsilon^{u}_t\psi_0 $ as $n$ tends to infinity.
\end{proposition}

\subsection{Controllability results}
Considerable efforts have been made to study the
controllability of~(\ref{eq:blse}). It is known (see \cite{turinici}) that exact controllability of (\ref{eq:blse}) in $H$ is hopeless in general. With the exception of some very particular examples where $\Omega$ is one dimensional, (\cite{camillo}), no description of the attainable set is known.  Therefore, we often consider the weaker notion of approximate controllability:
\begin{defn}
Let $(A,B)$ satisfy Assumption \ref{ass:ass}. The system $(A,B)$ is \emph{approximately controllable} if, for every $\psi_0,\psi_1$ in the unit Hilbert sphere, for every $\varepsilon>0$, there exists $u_{\varepsilon}:[0,T_{\varepsilon}]\rightarrow \mathbf{R}$ such that $\|\Upsilon^{u_\varepsilon}_{T_\varepsilon}\psi_0-\psi_1\|<\varepsilon$.
\end{defn}
Various methods have been used to give criterion of approximate controllability of system (\ref{eq:blse}). \cite{nersesyan}, \cite{beauchard-mirrahimi} and \cite{mirrahimi-continuous} rely on a Lyapunov approach.  \cite{Schrod2} adopt a more geometrical point of view, centered on the notion of non-degenerate (or non-resonant) transitions.
\begin{defn}
Let $(A,B)$ satisfy Assumption \ref{ass:ass}. A pair $(j,k)$ of integers is a  \emph{non-degenerate transition} of $(A,B)$ if (i)  $\langle \phi_j,B\phi_k \rangle \neq 0$ and (ii) for every $(l,m)$ in $\mathbf{N}^2$, $|\lambda_j-\lambda_k|=|\lambda_l-\lambda_m|$ implies $(j,k)=(l,m)$ or $\langle \phi_l, B \phi_m\rangle =0$ or $\{j,k\}\cap \{l,m\}=\emptyset$.
\end{defn}

 \begin{defn}
Let $(A,B)$ satisfy Assumption \ref{ass:ass}. A subset $S$ of $\mathbf{N}^2$ is a \emph{non-degenerate chain of connectedness} of $(A,B)$ if
(i) for every $(j,k)$ in $S$, $(j,k)$ is a non-degenerate transition of $(A,B)$ and (ii) for every $r_a,r_b$ in $\mathbf{N}$, there exists a finite sequence $r_a=r_0,r_1,\ldots,r_p=r_b$ in $\mathbf{N}$ such that, for every $j\leq p-1$, $(r_j,r_{j+1})$ belongs to $S$.   
\end{defn}

The following sufficient criterion for approximate controllability is the central result of \cite{Schrod2}.
\begin{prop}
 Let $(A,B)$ satisfy Assumption \ref{ass:ass}.  If $(A,B)$ admits a non-degenerate chain of connectedness, then $(A,B)$ is approximately controllable. 
\end{prop}
As proved by
\cite{genericity-mario-paolo} (see also \cite{genericity-mario-privat}), a system $(A,B)$ satisfying Assumption \ref{ass:ass} generically admits a non-degenerate chain of connectedness. Hence,  approximate controllability is a generic property for systems of the type
(\ref{eq:main}) (see also  \cite{nersesyan}).

In \cite{Schrod2} an extensive (and in some case implicit) 
use of averaging results has been made. The following result  is a generalization of the Rotating Wave Approximation to infinite dimensional systems and can be found in~\cite{periodic}.

\begin{prop}\label{PRO_contr_Schrod2}
Let $(A,B)$ satisfy Assumption \ref{ass:ass} and $(j,k)$ be a non-degenerate transition of $(A,B)$. Define  $T=2\pi/ |\lambda_j-\lambda_k|$ and $\mathcal{N}=\{(l,m)\in \mathbf{N}^2| \langle \phi_l,B\phi_m\rangle \neq 0 \mbox{ and } |\lambda_l-\lambda_m| \in (\mathbf{N}\setminus\{1\})|\lambda_j-\lambda_k| \mbox{ and } \{j,k\}\cap \{l,m\}\neq \emptyset \}$. Consider a  $T$-periodic function $u^{\ast}:\mathbf{R}\rightarrow \mathbf{R}$  satisfying $\int_0^T u^\ast(t) e^{\mathrm{i} (\lambda_j-\lambda_k) t} \mathrm{d}t \neq 0$ and $\int_0^T u^\ast(t) e^{\mathrm{i} (\lambda_l-\lambda_m) t} \mathrm{d}t = 0$ for every $(l,m)$ in $\mathcal N$ and let 
$$
T^\ast= \frac{\pi  T}{2 |b_{1,2}|  \left |\int_0^T \!\! u^{\ast}(\tau)e^{\mathrm{i}(\lambda_{1}-\lambda_{2}) \tau}  \mathrm{d}\tau \right |}.
$$ 
Then there exist a sequence $(T_n^\ast)_{n\in \mathbf{N}}$ such that $T_n^\ast \in (nT^\ast-T,nT^\ast+T)$ and
$|\langle  \phi_k, \Upsilon^{u^\ast/n}_{T^\ast_{n}}\phi_j\rangle |$ tends to one as $n$ tends to infinity.
\end{prop}

\subsection{Formulation of an optimal control problem}

Let $(A,B)$ satisfy Assumption \ref{ass:ass} and admit a non-degenerate chain of connectedness. For every $r>0$, for every $j,k$ in $\mathbf{N}$ and $\varepsilon>0$
we define
 $ \mathcal{A}_r^{\varepsilon}(j,k)$ as the set of  functions $u:[0,T_u]\rightarrow \mathbf{R}$ in $L^1([0,T_u])\cap L^r([0,T_u])$ such that $\|\Upsilon^{u}_{T_u}\phi_j-\phi_k\|<\varepsilon$.
We consider the quantity 
$$
\mathcal{C}_r(\phi_{j},\phi_{k})=\sup_{\varepsilon>0} \left ( \inf_{ u \in \mathcal{A}_r^{\varepsilon}(j,k)} \|u\|_{L^r(0,T_u)} \right).
$$  
This quantity is the infimum of the $L^{r}$-norm of a control achieving approximate controllability. It clearly satisfies the triangle inequality.
Next proposition states that $\mathcal{C}_{r}$ is a distance on the space of eigenlevels only when $r=1$. 
Its proof is given in Section~\ref{SEC_C_r}.
\begin{prop}\label{PRO_main_result}
$\mathcal{C}_1$ is a distance on the set $\{\phi_j,j\in \mathbf{N}\}$. 
For $r>1$, $\mathcal{C}_r$ is equal to zero on the set $\{\phi_j,j\in \mathbf{N}\}$.
\end{prop}

\subsection{Weakly-coupled systems}
\begin{definition}
Let $k$ be a positive number and let  $(A,B)$ satisfy Assumption \ref{ass:ass}.
Then $(A,B)$ is
\emph{$k$ weakly-coupled}
if for every $u_1 \in \mathbf{R}$, $D(|A+u_1B|^{k/2})=D(|A|^{k/2})$ and
 there exists
a constant $c_{(A,B)}$ such that, for every $\psi$ in $D(|A|^k)$, $ |\Re \langle |A|^k
\psi,B \psi \rangle |\leq c_{(A,B)} |\langle |A|^k \psi, \psi \rangle|$.
\end{definition}
The notion of weakly-coupled systems is closely related to the growth of the
 $|A|^{k/2}$-norm $\langle |A|^k \psi, \psi \rangle$.
 For $k=1$, this quantity is the expected value of the energy of the system.
\begin{proposition}{\cite[Proposition 2]{weakly-coupled}}\label{PRO_croissance_norme_A} Let  $(A,B)$ be $k$-weakly-coupled.  Then,
for every $\psi_{0} \in D(|A|^{k/2})$, $K>0$,
$T\geq 0$, and $u$ in $L^1([0,\infty))$  for which
$\|u\|_{L^1}< K$, one has
$\left\|\Upsilon^{u}_{T}(\psi_{0})\right\|_{k/2} \leq
e^{c(A,B) K} \| \psi_0 \|_{k/2}.$
\end{proposition}
For every $N$ in $\mathbf{N}$, we define $\mathcal{L}_N$ the linear space spanned by $\phi_1,\phi_2, \ldots,\phi_N$ and $\pi_N:H\rightarrow H$, the orthogonal projection onto $\mathcal{L}_N$:
$$
\pi_N(\psi)=\sum_{k=1}^N \langle \phi_k, \psi \rangle \phi_k.
$$
The compressions of order $N$ of $A$ and $B$ are the finite rank operators 
$A^{(N)}=\pi_N A_{{\upharpoonright \mathcal{L}_N}}$ and
$B^{(N)}=\pi_N B_{{\upharpoonright \mathcal{L}_N}}$. The Galerkin approximation of (\ref{eq:main}) at order $N$ is the infinite dimensional system
\begin{equation}\label{EQ_Galerkin}
 \frac{\mathrm{d}}{\mathrm{d}t}x=A^{(N)}x +u(t) B^{(N)} x.
\end{equation}
Since $\mathcal{L}_N$ is invariant by (\ref{EQ_Galerkin}), one may also consider (\ref{EQ_Galerkin}) as a finite-dimensional system, whose propagator is denoted by $X^u_{(N)}(t,s)$.

\begin{proposition}{\cite[Proposition 4]{weakly-coupled}}\label{prop:gga}
Let $k$ and $s$ be non-negative  numbers with
$0\leq s <k$. Let $(A,B)$ 
 be $k$ weakly-coupled
Assume that there
exists $d>0$, $0\leq r<k$ such that $\|B\psi \|\leq d \|\psi \|_{r/2}$ 
 for
every $\psi$ in $D(|A|^{r/2})$.
Then
for every $\varepsilon > 0 $, $K\geq 0$, $n\in \mathbf{N}$, and
$(\psi_j)_{1\leq j \leq n}$ in $D(|A|^{k/2})^n$
there exists $N \in \mathbf{N}$
such that
for every piecewise constant function $u$
$$
\|u\|_{L^{1}} < K \Rightarrow\| \Upsilon^{u}_{t}(\psi_{j}) -
X^{u}_{(N)}(t,0)\pi_{N} \psi_{j}\|_{s/2} < \varepsilon,
$$
for every $t \geq 0$ and $j=1,\ldots,n$.
\end{proposition}

\begin{remark}
 An interesting feature of Propositions \ref{PRO_croissance_norme_A} and \ref{prop:gga} is the fact that the bound  of the 
$|A|^{k/2}$ norm of the solution of (\ref{eq:main}) or the bound on the error between the infinite dimensional 
system and its finite dimensional approximation only depend on the $L^1$ norm of the control, not on the time. 
\end{remark}

\section{Proof of Proposition~\ref{PRO_main_result}}\label{SEC_C_r}

\subsection{Lower bounds for the $L^1$ norm}\label{SEC_L1norm_lower}
\begin{prop}\label{PRO_minoration_L1}
 Let $(A,B)$ satisfy Assumption \ref{ass:ass}. For every $j,k$ in $\mathbf{N}$  such that $B\phi_j\neq 0$, for every  locally integrable $u:[0,T]\rightarrow \mathbf{R}$, 
$$
\|u\|_{L^1(0,T)}\geq \frac{\big | |\langle \phi_j,\phi_k \rangle| - |\langle \phi_j,\Upsilon^u_T\phi_k \rangle |
\big |}{\|B\phi_j\|}.
$$ 
\end{prop}
\textbf{Proof:}
Let $j,k$ in $\mathbf{N}$. For every locally integrable $u:[0,T]\rightarrow \mathbf{R}$, 
define $y:u\mapsto e^{-A t} \Upsilon^u_t\phi_k$. For almost every $t$, $y$ is differentiable with respect to $t$ and
$\displaystyle{
\frac{\mathrm{d}y}{\mathrm{d}t}=u(t)e^{-At} B e^{At} y}$.
In particular,
$$
 \left | \left \langle \phi_j , \frac{\mathrm{d}y}{\mathrm{d}t} 
\right \rangle \right | \leq |u(t)| \|B \phi_j\|.
$$
This concludes the proof of Proposition \ref{PRO_minoration_L1}.  \hfill ~\qed

A consequence of Proposition \ref{PRO_minoration_L1} is that ${\mathcal C}_1(\phi_j,\phi_k)$ is bounded away from zero as soon as $j\neq k$.

\subsection{Upper bound for the $L^1$ norm}\label{SEC_L1norm_upper}

In order to give an upper bound for $\mathcal{C}_1(\phi_j,\phi_k)$ when $(j,k)$ is a non-degenerate transition of $(A,B)$, we come back to Proposition \ref{PRO_contr_Schrod2}.
In the case where $\mathcal N$ is finite, \cite{Schrod2} give an explicit construction of a piecewise constant $u^{\ast}$ with value in $[0,1]$ and satisfying the assumptions of Proposition \ref{PRO_contr_Schrod2} such that, for every $n$ in $\mathbf{N}$,  $$\left \| \frac{u^{\ast}}{n} \right \|_{L^1(0,T^\ast_n)}\leq \frac{5\pi}{4 |\langle \phi_j,B\phi_k\rangle |}.$$ 
More details about the choice of $u^{\ast}$ when $(A,B)$ is weakly-coupled are given by \cite[Section III-C]{FEPSacc}.
Let us just mention that the choice $u^{\ast}:t\mapsto \cos(|\lambda_j-\lambda|t)$ guarantees
 $$\left \| \frac{u^{\ast}}{n} \right \|_{L^1(0,T^\ast_n)}\leq \frac{2}{ |\langle \phi_j,B\phi_k\rangle |}$$ for every $n$ in $\mathbf{N}$.

In any case, $(A,B)$ being weakly-coupled or not, this guarantees that, for every $j,k$ in $\mathbf{N}$, $\mathcal{C}_1(\phi_j,\phi_k)<+\infty$.

The fact that $\mathcal{C}_1$ is symmetric is a consequence of the so-called  time reversibility of (\ref{eq:main}): the propagator associated with  $(-A,-B)$ is the adjoint of the  backward propagator  associated with $(A,B)$. 
This proves that $\mathcal{C}_1$ is a distance on the set $\{\phi_j,j\in \mathbf{N}\}$.

\subsection{$L^r$ norms with $r>1$}\label{SEC_Lpnorm}
\begin{prop}
Let $(A,B)$ satisfy Assumption \ref{ass:ass} and admit a non-degenerate chain of connectedness.  If $r > 1$, then, for every $j,k$ in $\mathbf{N}$, $\mathcal{C}_r(\phi_j,\phi_k)=0$. 
\end{prop}
\textbf{Proof:} It is enough to consider $(j,k)$ in a non-degenerate chain of connectedness of $(A,B)$. The result  is then a consequence of Proposition \ref{PRO_contr_Schrod2}, since, for every $n$ in $\mathbf{N}$,
\begin{eqnarray*}
\left \|\frac{u^{\ast}}{n} \right \|^r_{L^r(0,T^\ast_n)} &=&
\frac{1}{n^r} \int_0^{ T^\ast_n} \!\!\!\!\!\! |u^{\ast}(t)|^r \mathrm{d}t\\
&\leq & \frac{1}{n^r} \int_0^{ nT^\ast +T} \!\!\!\!\!\! |u^{\ast}(t)|^r \mathrm{d}t\\
&\leq & \frac{1}{n^r} \int_0^{ ( n \left \lceil \frac{T^\ast}{T} \right \rceil +1) T} \!\!\!\!\!\! |u^{\ast}(t)|^r \mathrm{d}t\\
&\leq & \frac{n}{n^r} \left (\frac{T^\ast}{T} +2  \right )\int_0^{T}\!\!\!\!\!\! |u^{\ast}(t)|^r \mathrm{d}t,
\end{eqnarray*}
which tends to zero as $n$ tends to infinity.\hfill ~\qed

\section{ROTATION OF A PLANAR MOLECULE}\label{SEC_example}
In this Section, we apply our results to the well studied example of the rotation of a planar molecule (see, for instance, \cite{salomon-HCN,noiesugny-CDC, Schrod2}).
\subsection{Presentation of the model}

We consider a linear molecule with fixed length and center of mass. We assume that the molecule is constrained to stay in a fixed plane and that its only degree of freedom is the rotation, in the plane, around its center of mass. The state of the system at time $t$ is described by a point $\theta \mapsto \psi(t,\theta)$ of $L^2(\Omega,\mathbf{C})$ where $\Omega=\mathbf{R}/2\pi \mathbf{Z}$ is the one dimensional torus.
The Schr\"odinger equation writes
\begin{equation}\label{EQ_Schrod_circle}
\mathrm{i} \frac{\partial \psi}{\partial t}(t,\theta)= -\Delta \psi(t,\theta) + u(t) \cos(\theta) \psi(t,\theta),
\end{equation}
where $\Delta$ is the Laplace-Beltrami operator on $\Omega$.
The  self-adjoint operator $-\Delta$ has purely discrete spectrum $\{k^2,k \in \mathbf{N}\}$. All its eigenvalues are double but zero which is simple. The eigenvalue zero is associated with the constant
functions. The eigenvalue $k^2$ for $k>0$ is associated with the two eigenfunctions $\theta \mapsto\frac{1}{\sqrt{\pi}} \cos(k \theta)$ and $\theta \mapsto \frac{1}{\sqrt{\pi}} \sin(k \theta)$. The Hilbert space $H=L^2(\Omega,\mathbf{C})$ splits in two subspaces $H_e$ and $H_o$, the spaces of even and odd functions of $H$ respectively. The spaces $H_e$ and $H_o$ are invariant under the dynamics of (\ref{EQ_Schrod_circle}), hence no global controllability is to be expected in $H$.

We focus on the space $H_o$. The restriction $A$ of $\mathrm{i} \Delta$ to $H_o$ is skew adjoint, with simple eigenvalues $(-\mathrm{i}k^2)_{k\in \mathbf{N}}$ associated with the eigenvectors 
$$
\left (\phi_k:\theta \mapsto \frac{1}{\sqrt{\pi}} \sin(k \theta) \right )_{k\in \mathbf{N}}.
$$ 
The restriction $B$ of $\psi\mapsto -\mathrm{i} \cos(\theta) \psi$ to $H_o$ is skew-adjoint and bounded. The pair $(A,B)$ satisfies Assumption \ref{ass:ass} and is weakly-coupled (see \cite[Section~III.C]{weakly-coupled}).

The Galerkin approximations of $A$ and $B$ of order $N$ are
$$
A^{(N)}= - \left ( \begin{array}{cccc}
                \mathrm{i} & 0 & \cdots & 0 \\
		0 	& 4 \mathrm{i} & \ddots & \vdots \\
		\vdots & \ddots & \ddots & 0\\
		0 & \cdots & 0 & N^2 \mathrm{i}
                \end{array}
\right )  \mbox{ ~~and }$$
$$
B^{(N)}= - \mathrm{i}\left ( \begin{array}{ccccc}
                0 & 1/2 & 0 & \cdots & 0 \\
		 1/2 & 0 & 1/2 & \ddots & \vdots \\
		0  & \ddots & 0  &\ddots  & 0\\
		\vdots & \ddots  & 1/2 & 0 & 1/2\\
		0 & \cdots & 0 & 1/2 & 0
                \end{array}
\right ).
$$

\subsection{Computation of $\mathcal{C}_1(\phi_1,\phi_2)$}
Our aim is to compute the minimal $L^1$ norm needed to approximately transfer the wave function from the first eigenspace to the second one. Precisely, we will prove
\begin{prop}\label{PRO_calcul_C_1}
$
{\mathcal{C}}_1(\phi_1,\phi_2)= \pi.
$
\end{prop}

\textbf{Proof:} The transition $(1,2)$ is non-degenerate. Proposition \ref{PRO_contr_Schrod2} applies with $\mathcal{N}=\emptyset$. For every $\eta$ in $(0,1)$, we define $u^{\eta}$, the $2\pi/3$ periodic function defined by
$$
\left \{
\begin{array}{lll}
 u^{\eta}(x)=1 & \mbox{ for } & 0<x<\eta\\
u^{\eta}(x)=0 & \mbox{ for } & \eta \leq x \leq 2\pi/3.
\end{array}
\right.
$$
Proposition \ref{PRO_contr_Schrod2} states that,
defining 
$$T^\ast=\frac{\pi^2}{\left | \sin \left ( \frac{3\eta}{2} \right ) \right |},$$
 there exists a sequence $(T_n^\ast)_{n\in \mathbf{N}}$ such 
that $T_n^\ast \in (nT^\ast-2\pi/3, nT^\ast +2\pi/3)$ and
$|\langle \phi_2, \Upsilon^{u^\eta/n}_{T^\ast_n} \phi_1 \rangle |$ tends to one as $n$ tends to infinity. 
  
One computes, for every $n$ in $\mathbf{N}$, 
$$
\left \| \frac{u^\eta}{n} \right \|_{L^1(0,T^\ast_n)} \leq \frac{\eta}{\left |\sin \left ( \frac{3\eta}{2} \right ) \right |} \frac{3\pi}{2}.
$$
This last quantity tends to $\pi$ as $\eta$ tends to zero. This proves that 
$
{\mathcal{C}}_1(\phi_1,\phi_2)\leq \pi.
$

For every locally integrable control $u$, we define $y_1:t\mapsto \langle \phi_1, \Upsilon^u_t\phi_1 \rangle$ and 
$y_2:t\mapsto \langle \phi_2, \Upsilon^u_t\phi_1 \rangle$. The function $t\mapsto |y_1(t)|^2$ is absolutely continuous, and, for almost every $t$,
\begin{eqnarray*}
 \frac{\mathrm{d}}{\mathrm{d}t}|y_1(t)|^2&=& 2\Re (y_1'(t) \bar{y}_1(t))\\
&=& u(t) \Re (y_2(t) \bar{y}_1(t)).
\end{eqnarray*}
Hence, 
$$
\left | \frac{\mathrm{d}}{\mathrm{d}t}|y_1(t)|^2 \right | \leq |y_2(t)||y_1(t)||u(t)|,
$$
or
$$
-|u(t)| \leq \frac{\frac{\mathrm{d}}{\mathrm{d}t}|y_1(t)|^2}{\sqrt{|y_1|^2} \sqrt{1-|y_1(t)|^2}} \leq |u(t)|
$$
Integrating between $0$ and $T$, one gets
$$
 2 \left |\arctan \left ( \sqrt{\frac{1}{|y_1(T)|^2}-1}\right ) \right |\leq \|u\|_{L^1(0,T)},
$$
and, provided $\|u\|_{L^1(0,T)} < \pi$,
$$
\sqrt{\frac{1}{|y_1(T)|^2}-1} \leq \tan \left (\frac{\|u\|_{L^1(0,T)}}{2} \right ).
$$
Finally,
$ \displaystyle{
|y_1(T)| \geq \cos \left ( \frac{\|u\|_{L^1(0,T)}}{2}\right )}$
and
$$
|y_2(T)| \leq \sqrt{1-|y_1(t)|^2} \leq \sin \left ( \frac{\|u\|_{L^1(0,T)}}{2}\right ).
$$
If $(u_n)_{n \in\mathbf{N}}$ is a sequence of locally integrable functions and 
$(T_n)_{n \in \mathbf{N}}$ is a sequence of positive numbers such that $|\langle \phi_2, \Upsilon^{u_n}_{T_n}\phi_1\rangle|$ 
tends to one, then $\liminf_n \|u_n\|_{L^1(0,T_n)}\geq \pi$, hence $\mathcal{C}_1(\phi_1,\phi_2)\geq \pi$.
This concludes the proof of Proposition \ref{PRO_calcul_C_1}.\hfill ~\qed

\section{CONCLUSION AND PERSPECTIVES}
\subsection{Conclusion}
We introduced an optimal control problem associated with a bilinear quantum system. For $p\geq 1$, we have given some estimates of the $L^p$ norm of the control needed to steer the system from an eigenstate of the free Hamiltonian to another. In particular, for generic bilinear quantum systems, it is possible to steer any eigenstate of  the free Hamiltonian to any neighborhood of any other eigenstate with arbitrary small $L^2$ norm of the control. 
\subsection{Future works}
The estimates given for the $L^1$ norm only depend on the control potential $B$ (and not on the eigenvalues of the free Hamiltonian $A$ as long as the transition stay non-degenerate). It is possible that $L^p$ costs, with $ p <1$, are physically more relevant. A new approach would be needed to study this case since 
our methods do not provide any information about the $L^p$  norm of the control needed to steer an energy  level to another in the case $p<1$.

\bibliography{biblio}

\begin{thebibliography}{15}
\providecommand{\natexlab}[1]{#1}
\providecommand{\url}[1]{\texttt{#1}}
\providecommand{\urlprefix}{URL }
\expandafter\ifx\csname urlstyle\endcsname\relax
  \providecommand{\doi}[1]{doi:\discretionary{}{}{}#1}\else
  \providecommand{\doi}{doi:\discretionary{}{}{}\begingroup
  \urlstyle{rm}\Url}\fi

\bibitem[{Beauchard and Laurent(2010)}]{camillo}
Beauchard, K. and Laurent, C. (2010).
\newblock Local controllability of 1{D} linear and nonlinear {S}chr{\"o}dinger
  equations with bilinear control.
\newblock \emph{J. Math. Pures Appl.}, 94(5), 520--554.

\bibitem[{Beauchard and Mirrahimi(2009)}]{beauchard-mirrahimi}
Beauchard, K. and Mirrahimi, M. (2009).
\newblock Practical stabilization of a quantum particle in a one-dimensional
  infinite square potential well.
\newblock \emph{SIAM J. Control Optim.}, 48(2), 1179--1205.
\newblock \doi{10.1137/070704204}.

\bibitem[{Boscain et~al.(2012)Boscain, Caponigro, Chambrion, and
  Sigalotti}]{Schrod2}
Boscain, U., Caponigro, M., Chambrion, T., and Sigalotti, M. (2012).
\newblock A weak spectral condition for the controllability of the bilinear
  {S}chr{\"o}dinger equation with application to the control of a rotating
  planar molecule.
\newblock \emph{Communications in Mathematical Physics}, 311(2), 423--455.

\bibitem[{Boscain et~al.(2009)Boscain, Chambrion, Mason, Sigalotti, and
  Sugny}]{noiesugny-CDC}
Boscain, U., Chambrion, T., Mason, P., Sigalotti, M., and Sugny, D. (2009).
\newblock Controllability of the rotation of a quantum planar molecule.
\newblock In \emph{Proceedings of the 48th IEEE Conference on Decision and
  Control}, 369--374.

\bibitem[{Boussa{\"{i}}d et~al.(2011{\natexlab{a}})Boussa{\"{i}}d, Caponigro,
  and Chambrion}]{FEPSacc}
Boussa{\"{i}}d, N., Caponigro, M., and Chambrion, T. (2011{\natexlab{a}}).
\newblock Periodic control laws for bilinear quantum systems with discrete
  spectrum.
\newblock Preprint arXiv 1111.4550, to be presented in ACC 2012.

\bibitem[{Boussa{\"{i}}d et~al.(2011{\natexlab{b}})Boussa{\"{i}}d, Caponigro,
  and Chambrion}]{weakly-coupled}
Boussa{\"{i}}d, N., Caponigro, M., and Chambrion, T. (2011{\natexlab{b}}).
\newblock Weakly-coupled systems in quantum control.
\newblock Preprint arXiv 1109.1900v1.

\bibitem[{Chambrion(2011)}]{periodic}
Chambrion, T. (2011).
\newblock Periodic excitations of bilinear quantum systems.
\newblock Preprint arXiv:1103.1130, to appear in \emph{Automatica}.

\bibitem[{Dahleh et~al.(1990)Dahleh, Peirce, and Rabitz}]{PhysRevA.42.1065}
Dahleh, M., Peirce, A.P., and Rabitz, H. (1990).
\newblock Optimal control of uncertain quantum systems.
\newblock \emph{Phys. Rev. A}, 42, 1065--1079.
\newblock \doi{10.1103/PhysRevA.42.1065}.
\newblock \urlprefix\url{http://link.aps.org/doi/10.1103/PhysRevA.42.1065}.

\bibitem[{Grivopoulos and Bamieh(2008)}]{symeon}
Grivopoulos, S. and Bamieh, B. (2008).
\newblock Optimal population transfers in a quantum system for large transfer
  time.
\newblock \emph{IEEE Trans. on Autom. Control}, 53(4), 980--992.

\bibitem[{Mason and Sigalotti(2010)}]{genericity-mario-paolo}
Mason, P. and Sigalotti, M. (2010).
\newblock Generic controllability properties for the bilinear {S}chr{\"o}dinger
  equation.
\newblock \emph{Communications in Partial Differential Equations}, 35,
  685--706.

\bibitem[{Mirrahimi(2009)}]{mirrahimi-continuous}
Mirrahimi, M. (2009).
\newblock Lyapunov control of a quantum particle in a decaying potential.
\newblock \emph{Ann. Inst. H. Poincar{\'e} Anal. Non Lin{\'e}aire}, 26(5),
  1743--1765.
\newblock \doi{10.1016/j.anihpc.2008.09.006}.

\bibitem[{Nersesyan(2010)}]{nersesyan}
Nersesyan, V. (2010).
\newblock Global approximate controllability for {S}chr{\"o}dinger equation in
  higher {S}obolev norms and applications.
\newblock \emph{Ann. Inst. H. Poincar{\'e} Anal. Non Lin{\'e}aire}, 27(3),
  901--915.
\newblock \doi{10.1016/j.anihpc.2010.01.004}.

\bibitem[{Privat and Sigalotti(2010)}]{genericity-mario-privat}
Privat, Y. and Sigalotti, M. (2010).
\newblock The squares of the {L}aplacian--{D}irichlet eigenfunctions are
  generically linearly independent.
\newblock \emph{ESAIM: COCV}, 16, 794--807.

\bibitem[{Salomon and Turinici(2005)}]{salomon-HCN}
Salomon, J. and Turinici, G. (2005).
\newblock Control of molecular orientation and alignment by monotonic schemes.
\newblock In \emph{Proceedings of the 24-th IASTED International Conference on
  modelling, identification and control}, 64--68.

\bibitem[{Turinici(2000)}]{turinici}
Turinici, G. (2000).
\newblock On the controllability of bilinear quantum systems.
\newblock In M.~Defranceschi and C.~{Le Bris} (eds.), \emph{Mathematical models
  and methods for ab initio Quantum Chemistry}, volume~74 of \emph{Lecture
  Notes in Chemistry}. Springer.

\end{thebibliography}

\end{document}